\documentclass[psamsfonts]{article}

    \usepackage{latexsym}
    \usepackage{amsfonts}
    \usepackage{amsmath}
    \usepackage{amssymb}
    \usepackage{amscd}
    \usepackage{amstext}
    \usepackage{amsthm}
    \usepackage{color}

    \title{Tracial Rokhlin property for
    automorphisms on simple $A{\mathbb T}$-algebras}

    \author{Huaxin Lin\\
    Department of Mathematics\\
    University of Oregon\\
    Eugene, Oregon 97403-1222\\
    Hiroyuki Osaka\\
    Department of Mathematical Sciences\\
    Ritsumeikan University\\
    Kusatsu, Shiga, 525-8577, Japan
    }
    \date{}

\begin{document}
    \maketitle


    \newcommand{\CA}{$C^*$-algebra}
    \newcommand{\SCA}{$C^*$-subalgebra}
\newcommand{\aue}{approximate unitary equivalence}
    \newcommand{\ayue}{approximately unitarily equivalent}
    \newcommand{\mops}{mutually orthogonal projections}
    \newcommand{\hm}{homomorphism}
    \newcommand{\pisca}{purely infinite simple \CA}
    \newcommand{\andeqn}{\,\,\,\,\,\, {\rm and} \,\,\,\,\,\,}
    \newcommand{\QED}{\rule{1.5mm}{3mm}}
    \newcommand{\morp}{contractive completely
    positive linear map}
    \newcommand{\asmorp}{asymptotic morphism}
    \newcommand{\arrow}{\rightarrow}
    \newcommand{\tdsum}{\widetilde{\oplus}}
    \newcommand{\pa}{\|}  
    \newcommand{\ep}{\varepsilon}
    \newcommand{\id}{{\rm id}}
    \newcommand{\aueeps}[1]{\stackrel{#1}{\sim}}
    \newcommand{\aeps}[1]{\stackrel{#1}{\approx}}
    \newcommand{\dt}{\delta}
    \newcommand{\yu}{\fang}
    \newcommand{\ca}{{\cal C}_1}
\newcommand{\Ad}{{\rm ad}}
    \newcommand{\tr}{{\rm TR}}
    \newcommand{\N}{{\bf N}}
    \newcommand{\C}{{\bf C}}
    \newcommand{\Aut}{{\rm Aut}}

\newcommand{\cg}{\textcolor[rgb]{0.00,0.47,0.45}}

\newcommand{\cp}{\textcolor[rgb]{1.00,0.00,1.00}}

\newcommand{\T}{{\mathbb T}}
\newcommand{\R}{{\mathbb R}}
    \newtheorem{thm}{Theorem}[section]
    \newtheorem{Lem}[thm]{Lemma}
    \newtheorem{Prop}[thm]{Proposition}
    \newtheorem{Dfn}[thm]{Definition}
    \newtheorem{Cor}[thm]{Corollary}
    \newtheorem{Ex}[thm]{Example}
    \newtheorem{Pro}[thm]{Problem}
    \newtheorem{Remark}[thm]{Remark}
    \newtheorem{NN}[thm]{}
    \renewcommand{\theequation}{e\,\arabic{section}.\arabic{equation}}
    \newcommand{\rforal}{{\rm\,\,\,for\,\,\,all\,\,\,}}
\newcommand{\Z}{{\mathbb Z}}
\newcommand{\Q}{{\mathbb Q}}

    \newcommand{\Ik}{ {\cal I}^{(k)}}
    \newcommand{\Iz}{{\cal I}^{(0)}}
    \newcommand{\Ii}{{\cal I}^{(1)}}
    \newcommand{\Ip}{{\cal I}^{(2)}}
\newcommand{\af}{\alpha}
\newcommand{\bt}{\beta}
    \renewcommand{\normalbaselines}{\baselineskip20pt
    \lieskip3pt \lineskiplimit3pt}

    \newcommand{\mapright}[1]{\smash{\mathop{
    \hbox to 1cm{\rightarrowfill}}\limits^{#1}}}

    \newcommand{\mapdown}[1]{\Big\downarrow
    \llap{$\vcenter{\hbox{$\scriptstyle#1\,$}}$ }}

\begin{abstract}
Let $A$ be a unital simple $A\T$-algebra of real rank zero.
Given an
isomorphism $\gamma_1: K_1(A)\to K_1(A),$ we show that there is an automorphism
$\af: A\to A$ such that $\af_{*1}=\gamma_1$ which has the tracial Rokhlin property.
Consequently, the crossed product $A\rtimes_{\af}\Z$ is a simple unital AH-algebra with real rank zero.
We also show that automorphism   with Rokhlin property can be constructed
from minimal homeomorphisms on a connected compact metric space.

\end{abstract}

\section{Introduction}
Alan Connes introduced the Rokhlin property in ergodic theory to operator algebras (\cite{Con}).
Several versions of the Rokhlin property for automorphisms on \CA s have been studied (for example
\cite{HO}, \cite{Ro}, \cite{Ks1}, \cite{Nk}, \cite{Iz}, \cite{Iz2}  and \cite{P1}, 
to name a few).
Given a unital \CA\, $A$ and an automorphism $\af$ on $A,$ one may view the pair $(A,\af)$ as
a non-commutative dynamical system. To study its dynamical structure, it is natural 
to introduce the notion of Rokhlin property. Let $A$ be a unital simple $A\T$-algebra (a \CA\, which
is an inductive limit of those \CA s that are finite direct sums of continuous functions
on the circle $\T$) and let $\af$ be an automorphism on $A.$ 
In several cases, Kishimoto showed that
if $\af$ is approximately inner (or homotopic to the identity) and has a Rokhlin type property,
then the crossed product $A\rtimes_{\af}\Z$ is again
a unital simple $A\T$ of real rank zero (\cite{Ks1}, \cite{Ks2}, \cite{Ks3}).
It is proved that if $\af$ induces the identity on $K_0(A)$ (or a ``dense"
subgroup of $K_0(A)$) and $\af$ has so-called tracial cyclic Rokhlin property then indeed the crossed
product is an $AH$-algebra of real rank zero (\cite{LO} and \cite{Lngpots}).
A natural question is when automorphisms
have certain Rokhlin property. In \cite{OP2}, it is shown that if $A$ is a unital separable simple \CA\, with tracial rank zero and with a unique tracial state and if $\af$ is an automorphism such that $A\rtimes_{\af}\Z$ has a unqiue trace,
then $\af$ has the tracial Rokhlin property.
More recently, N.C. Phillips \cite{Phn} showed that, for any unital simple
separable \CA\, $A$ with tracial rank zero, there is a dense $G$-$\dt$ set of approximately inner automorphisms such that every automorphism in the the set has the tracial Rokhlin property.

 Let $A$ be a unital simple $A\T$-algebra of real
rank zero. Suppose that $\gamma_1: K_1(A)\to K_1(A)$ is an automorphism.
In this note, we present an automorphism $\af$ on $A$ with the tracial cyclic Rokhlin property
such that $\af_{*1}=\gamma_1.$
The automorphism $\af$ that we present in this note also has the property that $\af_{*0}={\rm id}_{K_0(A)}.$
When $\af$ has the tracial cyclic Rokhlin property,
then $\af$ must fix a large subgroup of $K_0(A).$
In fact, it is shown in \cite{Lngpots}  that, at least in the case that $A$ has a unique tracial state, if $\af$ has tracial cyclic Rokhlin property,
then $\af^2$ fixes a subgroup $G\subset K_0(A)$
so that $\rho(G)$ is dense in $Aff(T(A)),$
where $T(A)$ is the tracial state space of $A,$ $Aff(T(A))$
is the space of all real affine continuous functions on $T(A)$ and
$\rho: K_0(A)\to Aff(T(A))$ is the positive
\hm\, induced by the evaluation
$\rho([p])(\tau)=\tau(p)$ for projections $p\in A.$

Let $X$ be a connected compact metric space and
$\rho_X: K_0(C(X))\to \Z$ be the dimension map.
Denote by ${\rm ker}\rho_X$ the kernel of the dimension map.
Given a such  $X$ and a countable dense subgroup
$D\subset \Q,$  there is a standard way to construct a unital simple
\CA \, $A_X$ with tracial rank zero such that
$K_0(A_X)=D\oplus {\rm ker}\rho_X$ and $K_1(A_X)=K_1(C(X))$  as well as there is a unital  embedding $j: C(X)\to A_X$ 
so that $(j_{*0})|_{{\rm ker}\rho_X}={\rm id}|_{{\rm ker}\rho_X}$ and 
$j_{*1}={\rm id}_{K_1(X)}.$
This could be viewed as a version of the non-commutative space associated with $X.$
Suppose that $\psi: X\to X$ is a minimal homeomorphism on $X.$
We show that one can construct an automorphism $\af$ on $A_X$ associated with $\psi$ such that
$\af$ has the tracial cyclic Rokhlin property and such that
$\af_{*0}|_D={\rm id}_D,$
$\af_{*0}|_{{\rm ker}\rho_X}=\psi_{*0}|_{{\rm ker}\rho_X}$ and
$\af_{*1}=\psi_{*1}.$ This may be viewed as a non-commutative version of the minimal
action associated with $\psi.$ We also show that somewhat general construction
can also be made. It appears that automorphisms with the tracial cyclic Rokhlin property
occur quite often. While we believe that many other types of construction of automorphisms
with the Rokhlin property may be possible and perhaps not necessarily difficult, we think these construction
in this note shed some  light on how commutative Rokhlin tower lemma appears naturally
in the study of non-commutative dynamical systems such as $(A,\af).$

\vspace{0.2in}

{\bf Aknowledgement}
The first named author was partially supported by a NSF grant.
The second author was partially supported by Open Research Center
Project for Private Universities:maching fund from MEXT,2004-2008.
Much of the ground work
of this reseach was done when both authors were visiting East China Normal University
in the summer 2004. They would like to acknowledge the support from Shanghai Priority Academic
Disciplines and from Department of Mathematics of East China Normal University.

\section{Preliminaries}

The following conventions will be used in this paper.
\vspace{0.2in}

(1) Let $A$ be a stably finite \CA. Denote by $T(A)$ the tracial
state space of $A.$ Denote by $Aff(T(A))$ the normed space of all
real affine continuous functions on $T(A).$

\vspace{0.1in}

 (2)
Denote by $\rho_A: K_0(A)\to Aff(T(A))$ the positive \hm\, induced
by $\rho_A([p])(\tau)=\tau\otimes Tr(p)$ for any projection in $M_k(A)$
($k=1,2,...,$), where $Tr$ is the standard trace on $M_k.$

\vspace{0.1in}

 (3) Let $X$ be a connected compact metric space.
Denote by $\rho_X: K_0(C(X))\to \Z$ the positive \hm\,
$\rho_{C(X)}$ defined in (2). It is the dimension function from
$K_0(C(X))$ to $\Z.$

\vspace{0.1in}

(4) Let $X$ be a compact metric space, let $F$ be a subset of $X$
and let $\ep>0.$ Put
$$
F_{\ep}=\{x\in X: {\rm dist}(x,F)<\ep\}.
$$

\vspace{0.1in}

(5) Let $A$ be a \CA\, and let $p$ and $q$ be projections in $A.$
We say $p$ is equivalent to $q$ if there is a $v\in A$ such that
$v^*v=p$ and $vv^*=q.$

\vspace{0.1in}

(6) Let $A=\lim_{n\to\infty}(A_n, \phi_{n})$ be an inductive limit
of \CA s. Here $\phi_n$ is a \hm\, from $A_n$ into $A_{n+1}.$ We
will use $\phi_{n, \infty}: A_n\to A$ for the \hm\, induced by the
inductive limit system. We also use $\phi_{n, m}$ for
$\phi_{m-1}\circ \phi_{m-2}\circ\cdots \circ \phi_n,$ if $m>n.$
Suppose that each $A_n=\oplus_{i=1}^{r(n)}B_{n,i}.$ When $B_{n,i}$
are understood, a partial map $\phi_n^{(i,j)}$ of $\phi_{n}$ is
the \hm\, from  $B_{n,i}$ to $B_{n+1,j}$ given by $\phi_n.$

\vspace{0.1in}

 (7) Let $A$ and $B$ be two \CA s and let
$\phi,\psi: A\to B$ be two maps. Suppose that ${\cal G}\subset A$
is a subset. We write
$$
\psi\approx_{\ep} \phi\,\,\,{\rm on}\,\,\,{\cal G},
$$
if
$$
\|\phi(a)-\psi(a)\|<\ep\rforal a\in {\cal G}.
$$

\vspace{0.2in}

The following Rokhlin property was introduced in \cite{OP1}.

\begin{Dfn}\label{Drok0}
{\rm
Let $A$ be a unital \CA\, and let $\af$ be an automorphism on $A.$
We say that $\af$ has the tracial Rokhlin property if, for any $\ep>0,$
any finite subset ${\cal F}\subset A,$ any positive integer $n>0,$ and any $a\in A_+\setminus\{0\},$
there are mutually orthogonal projections $e_1,e_2,...,e_{n}$ such that

(1) $\|e_ix-xe_i\|<\ep$ for all $x\in {\cal F},$
\ $(0 \leq i \leq n - 1)$

(2) $\|\af(e_i)-e_{i+1}\|<\ep,$ $i=1,2,...,n-1$ and

(3) $1-\sum_{i=0}^{n-1}e_i$ is equivalent to a projection in ${\overline{aAa}}.$
}
\end{Dfn}

A stronger version of the Rokhlin property below was given in \cite{LO}.

\begin{Dfn}\label{Drok1}
{\rm Let $A$ be a unital \CA\, and let $\af$ be an automorphism on $A.$
We say that $\af$ has the tracial {\it cyclic} Rokhlin property if, for any $\ep>0,$
any finite subset ${\cal F}\subset A,$ any positive integer $n>0,$ and any $a\in A_+\setminus\{0\},$
there are mutually orthogonal projections $e_0,e_1,e_2,...,e_{n-1}$ (with $e_n=e_0$) such that

(1) $\|e_ix-xe_i\|<\ep$ for all $x\in {\cal F},$
\ $(0 \leq i \leq n - 1)$

(2) $\|\af(e_i)-e_{i+1}\|<\ep,$ $i=0,1,...,n-1$ and

(3) $1-\sum_{i=0}^{n-1}e_i$ is equivalent to a projection in ${\overline{aAa}}.$
}
\end{Dfn}

\begin{Remark}
{\rm 
Note that, in Definition \ref{Drok1}, 
one requires that $\|\af(e_{n-1})-e_0\|<\ep.$ 
This is not required in Definition \ref{Drok0}.
If $A$ is assumed to satisfy the so-called fundamental comparison
property, i.e., $\tau(p)>\tau(q)$ for all $\tau\in T(A)$ implies that $q$ is equivalent to a
projection $e\le p$ for any pair of projections, then condition (3) above can be replaced
by $\tau(1-\sum_{i=0}^{n-1}e_i)<\ep$ for all $\tau\in T(A)$ and the reference to the
element $a$ can be removed.  We note that if  $A$ has tracial rank zero then $A$ has the fundamental comparison
property.
}
\end{Remark}

If $A$ has tracial rank zero, then $A$ has real rank zero, stable
rank one, and $\rho_A(K_0(A))$ is dense in $Aff(T(A)).$ One of the
important consequences that a given automorphism satisfies the
tracial cyclic  Rokhlin property is the following result:

\begin{thm}\label{gpots}{\rm (see  3.4 and 4.5 of \cite{Lngpots})}

If $A$ is a unital separable amenable simple \CA\, with tracial rank zero
satisfying the UCT
and $\af$ is an automorphism on $A$ with the tracial (cyclic)
Rokhlin property such that $(\af)_{*0}={\rm id}_G$ for some
subgroup $G\subset K_0(A)$ for which $\rho_A(G)$ is dense in
$\rho_A(K_0(A)),$ then $A\rtimes_{\af}\Z$ has tracial rank zero.

Consequently, $A\rtimes_{\af}\Z$ is an AH-algebra.
\end{thm}

\section{Automorphisms on simple $A\T$-algebras}
The purpose of this section is to present Theorem \ref{MT1}.

We start with the following construction which is certainly familiar to experts.

\begin{Lem}\label{smallpro}
Let $A$ be a unital simple infinite dimensional AF-algebra. Then,
there exists a sequence of finite dimensional \SCA s
$F_n=M_{k(n)}\oplus B_n$ and monomorphisms $\phi_n: F_n\to
F_{n+1}$ satisfying the following:

{\rm (1)} Let $d_n$ be the identity of $M_{k(n)}$. Then
$$
\lim_{n\to\infty}\sup_{\tau\in T(A)}\{\tau(\phi_{n,
\infty}(d_n))\}=0;
$$

{\rm (2)} The partial map from $M_{k(n)}$ to $M_{k(n+1)}$ has
multiplicity at least 2.

{\rm (3)} $A=\lim_{n\to\infty}(F_n, \phi_n).$

\end{Lem}

\begin{proof}
We first write $A=\lim_{n\to\infty}(A_n, \psi_n),$ where each
$A_n$ is a finite dimensional \CA\, and $\psi_n$ is a
monomorphism.

Write $A_1=M_{n(1)}\oplus B_1'.$ Since $A$ is assumed to be
simple, we may also assume that the partial map $\psi_n^{1,1}$
from $M_{n(1)}$ to a simple summand $M_{n(2)}$ of $A_2$ has
multiplicity $m(1)$ which is at least $4.$ Define
$F_1=M_{n(1)}\oplus M_{n(1)}\oplus\cdots \oplus M_{n(1)}\oplus
B_1',$ where $M_{n(1)}$ repeats $[m(1)/2]$ (integer part of
$m(1)/2$) times. Denote $B_1=M_{n(1)}\oplus\cdots \oplus
M_{n(1)}\oplus B_1',$ where $M_{n(1)}$ repeats $[m(1)/2]-1$ times.
So $F_1=M_{n(1)}\oplus B_1.$
Note there are monomorphism $f_1: A_1\to F_1$ and monomorphism
$\psi_1': F_1\to A_2$ such that $\psi_1=\psi_1'\circ f_1$ and the
multiplicity of the partial map of $\psi_1'$ from $M_{n(1)}$ to
$M_{n(2)}$ is at least 2.  Let $d_1$ be the identity of $M_{n(1)}$
in $F_1.$ Note that
$$
\tau(\psi_{2,\infty}\circ \psi_1'(d_1))\le 1/2
$$
for all $\tau\in T(A).$

We have $A_2=M_{n(2)}\oplus B_2',$ where $B_2'$ is a finite
dimensional \CA.

We may also assume (by replacing $A_3$ by some $A_n$) that the
partial map of $\psi_2$ from $M_{n(2)}$ to a simple summand
$M_{n(3)}$ of $A_3$ has multiplicity $m(2)$ which is at least
least $8.$ Define $F_2=M_{n(2)}\oplus M_{n(2)}\oplus\cdots
M_{n(2)}\oplus B_2',$ where $M_{n(2)}$ repeats $[m(2)/2]$ times.
Define $B_2=M_{n(2)}\oplus\cdots \oplus M_{n(2)}\oplus B_2',$ where
$M_{n(2)}$ repeats $[m(2)/2]-1$ times. So $F_2=M_{n(2)}\oplus B_2.$
There exists a monomorphism
$f_2: A_2\to F_2$ and there exists a monomorphism $\psi_2': F_2\to
A_3$ such that
$$
\psi_2=\psi_2'\circ f_2.
$$
Moreover, the multiplicity of the partial map of $\psi_2'$ from
$M_{n(2)}$ to $M_{n(3)}$ is at least $2.$

 Define $\phi_1: F_1\to F_2$ by $\phi_1=f_2\circ \psi_1'.$ Let
$d_2$ be the identity of the first summand $M_{n(2)}$ in $F_2.$
Then
$$
\tau(\psi_{3, \infty}\circ \psi_2'(d_2))<1/4
$$
for all $\tau\in T(A).$

Write $A_3=M_{n(3)}\oplus B_3',$ where $B_3'$ is a finite
dimensional \SCA. We may assume that the partial map of $\psi_3$
from $M_{n(3)}$ to a simple summand $M_{n(4)}$ of $A_4$ has
multiplicity $m(3)$ at least $16.$ Define $F_3=M_{n(3)}\oplus
M_{n(3)}\oplus\cdots \oplus M_{n(3)}\oplus B_3',$ where $M_{n(3)}$
repeats $[m(3)/2]$ times and $B_3=M_{n(3)}\oplus\cdots \oplus
M_{n(3)}\oplus B_3',$ where $M_{n(3)}$ repeats $[m(3)/2]-1$ times.
So $F_3=M_{n(3)}\oplus B_3.$

There are monomorphisms $f_3: A_3\to F_3$ and  $\psi_3': F_3\to
A_4$ such that $\psi_3=\psi_3'\circ f_3.$ Moreover, the partial
map of $\psi_3$ from $M_{n(3)}$ (the first summand) to $M_{n(4)}$
has multiplicity at least $2.$ Define $\phi_2: F_2\to F_3$ by
$\phi_2=f_3\circ \psi_2'.$ Let $d_3$ be the identity of of the
first summand $M_{n(3)}$ in $F_3.$ Then
$$
\tau(\psi_{4, \infty}\circ\psi_3'(d_3))<1/8.
$$

We continue this construction. We obtain $F_k=M_{n(k)}\oplus B_k,$
where $B_k$ is a finite dimensional \SCA, and we obtain monomorphisms $f_k:
A_k\to F_k,$ $\psi_k': F_k\to A_{k+1}$ and $\phi_k : F_k\to
F_{k+1}$ such that
\begin{eqnarray}\label{eadd1}
 \psi_k=\psi_k'\circ f_k \andeqn
\phi_k=f_{k+1}\circ \psi_k'
\end{eqnarray}
 and such that the partial
map of $\phi_k$ from $M_{n(k)}$ to $M_{n(k+1)}$ has multiplicity
at least $2.$ Moreover,
\begin{eqnarray}\label{eadd2}
\tau(\psi_{k+1,\infty}\circ \psi_k'(d_k))<1/2^k
\end{eqnarray}
for all $\tau\in T(A),$ where $d_k$ is the identity of $M_{n(k)}$
in $F_k.$

By (\ref{eadd1}), $A=\lim_{n\to\infty}(F_k, \phi_k).$ By
(\ref{eadd2}),
$$
\tau(\phi_{k, \infty}(d_k))<1/2^k.
$$
The lemma follows.

\end{proof}

\begin{Lem}\label{LLambda}
Let $X=S^1\vee S^1\vee\cdots \vee S^1$ be identified with $m$
copies of unit circle with a common point $1.$  For any integer $n$ and
finite subset ${\cal F}\subset C(X),$  there exists an integer
$N=\Lambda(n, {\cal F})$ satisfying the following:

For any unital \CA\ $B$
with real rank zero and any
\hm s  $\phi_1,\phi_2: C(X)\to B$
with $(\phi_1)_{*1}=(\phi_2)_{*1},$ there
is a unitary $u\in M_{mN+1}(B)$ such that
$$
u^*({\rm diag}(\phi_1(f), f(\xi_1),f(\xi_2),...,f(\xi_{mN}))u
\approx_{1/2^n}
$$
$$
 {\rm diag}(\phi_2(f), f(\xi_1),f(\xi_2),...,
f(\xi_{mN})),
$$
where
each of $\{\xi_1,\xi_2,...,\xi_N\},$
$\{\xi_{N+1},\xi_{N+2},...,\xi_{2N}\},$ $...,$
$\{\xi_{(m-1)N+1},\xi_{(m-1)N+2},...,\xi_{mN}\},$
 divides the unit
circles evenly.
\end{Lem}

\begin{proof}
This follows from Theorem 1.1 in \cite{GL} and its remark.

\end{proof}

\begin{Dfn}\label{Standard}
{\rm Let $X=S^1\vee S^1\vee\cdots \vee S^1$  and $Y=S^1\vee
S^1\vee\cdots \vee S^1$ be identified with $m$ and $r$ copies of
unit circle with a common point $1,$ respectively. Let $\pi_i:
Y\to S^1$ (the $i$th copy of $Y$) be defined as follows:
$$
\pi_i((\xi_1,\xi_2,...,\xi_i, ...,\xi_r))=\xi_i.
$$


A continuous map $s: X\to Y$ is said be {\it standard} if
 $\pi_i\circ s$ on each $S^1$
is described as $z\to z^k$ for some  integer $k,$ where $z$ is
the identity map on the unit circle.
Note that if $k=0,$ $z\to z^0$ is the map to the common point $1.$

Suppose that $\kappa: \Z^r\to \Z^m$ is the \hm\, associated with a
$m\times r$ matrix $(c_{i,j})$ with integer entries. Then we say
that $s: X\to Y$ is the {\it standard map associated with
$\kappa$} if $\pi_j\circ s$ on the $i$-th copy of $S^1$ can be
described by $z\to z^{c_{i,j}}.$}

\end{Dfn}

Note that, for a standard map $s: X\to Y,$ $s$ maps the common
point $1$ of $X$ to the common point $1$ of $Y.$

It should be also noted that if $Z=S^1\vee S^1\vee \cdots \vee
S^1$ is identified with $k$ copies of the unit circle with a
common point $1$ and $s': Y\to Z$ is  a standard map, then
$s'\circ s$ is a standard map from $X$ to $Z.$

\begin{Lem}\label{smallpro2}
Let $A$ be a unital simple A$\T$-algebra with real rank zero.
Suppose that $K_1(A)=\lim_{n\to\infty}(G_n, \kappa_n)\not=\{0\},$
where $G_n$ is direct sum of $\gamma(n)$ copies of $\Z$ and
$\kappa_n: G_n\to G_{n+1}$ is an \hm.  Then
$$
A=\lim_{n\to\infty}(A_n, \phi_n)
$$
such that

 {\rm (1)} $A_n=C(X_n)\otimes M_{k(n)}\oplus B_n,$ where $B_n$ is a
 finite dimensional \SCA s,

 {\rm (2)} $X_n=S^1\vee S^1\vee \cdots \vee S^1,$ where
 $S^1$ repeats $\gamma(n)$ times,

 {\rm (3)}
 $$
\lim_{n\to\infty}\sup_{\tau\in T(A)}\{\tau(\phi_{n,
\infty}(d_n))\}=0,$$ where $d_n$ is the identity of $C(X_n)\otimes
M_{k(n)},$

 {\rm (4)} each $\phi_n$ is
 injective.

Moreover, if $\{X_n\}$ is fixed and a sequence of positive
integers $\{a_n\}$ is given, we may require that

{\rm (5)} the partial map of $\phi_n$ from $C(X_n)\otimes
M_{k(n)}$ to $C(X_{n+1})\otimes M_{k(n+1)}$ has the form
$$
{\phi_n(f)={\rm diag}(f\circ s_n, f(\zeta(n,1)),f(\zeta(n,2)),...,
f(\zeta(n,t(n)))),}
$$
where $s_n$ is a standard map associated with $\kappa_n$ and
$\{\zeta(n,1),\zeta(n,2),...,\zeta(n,t(n))\}$ contains a subset that divides the $j$-th
circle into $b(n,j)$ even arcs with $b(n,j)\ge a_n;$

 {\rm (6)} each partial map $\phi_n^{(1,i)}$ of $\phi_n$ from
 $C(X_n)\otimes M_{k(n)}$ to a summand of $B_{n+1}$ has the form
 $$
{\phi_n^{(1,i)}(f)=\oplus_{j=1}^{\gamma(n)}{\rm diag} ( f(\xi(j,1)),
f(\xi(j,2)),...,f(\xi(j,k_j))),}
$$
where  $\{\xi(j,1),...,\xi(j,k_j)\}$ is $k_j$
points on the  $j$-th circle of $X_n.$

 \end{Lem}
\begin{proof}
Note that $A$ is an infinite dimensional simple \CA.
Let $B$ be a unital simple separable AF-algebra such that
$$
(K_0(B), K_0(B)_+, [1_B])=(K_0(A), K_0(A)_+, [1_A]).
$$
By Lemma \ref{smallpro}, we may assume that
$B=\lim_{n\to\infty}(M_{k(n)}\oplus B_n,\phi_n),$ where $B_n$ is a
finite dimensional \CA\, and (1) and (2) in Lemma \ref{smallpro} hold.
By passing to a subsequence, one may assume that the partial map
of $\phi_n$ from $M_{k(n)}$ to $M_{k(n+1)}$ has multiplicity at
least $m(n)\ge \gamma(n)a_n-1.$

Put $X_n=S^1\vee S^1\vee\cdots
\vee S^1,$ where $S^1$ repeats $\gamma(n)$ times. Define
$A_n=C(X_n)\otimes M_{k(n)}\oplus B_n.$ So
$K_1(A_n)=G_n=\Z^{\gamma(n)},$ $n=1,2,....$ Define $s_n: X_n\to
X_{n+1}$ to be the standard map associated with $\kappa_n$ as
defined in Definition \ref{Standard}. Define $\psi_n^{(1,0)}$ from
$C(X_n)\otimes M_{k(n)}$ to $C(X_{n+1})\otimes M_{k(n+1)}$ as
follows:
$$
\psi_n^{(1,0)}(f)={\rm diag}(f\circ s_n,
f(\zeta(n,1)),f(\zeta(n,2)),..., f(\zeta(n,t(n))))
$$
where $t(n)=m(n)-1$ and
$\{\zeta(n,1),\zeta(n,2),...,\zeta(n,t(n))\}$ contains a subset
that divides the $j$-th circle into $b(n,j)$ even arcs with
$b(n,j)\ge a_n$ for each $j.$

By identifying $M_{k(n+1)}$ with ${\mathbb C}\otimes M_{k(n+1)},$
we also assume that $\psi_n^{(1,0)}(1_{C(X_n)\otimes
M_{k(n)}})=\phi_n^{(1,0)}(1_{M_{k(n)}}),$ where $\phi_n^{(1,0)}$
is the partial map of $\phi_n$ from $M_{k(n)}$ to $M_{k(n+1)}.$
Let $m(n,1,i)$ be the multiplicity of the partial map of $\phi_n$
from $M_{k(n)}$ to the $i$-th summand of $B_{n+1}.$ We may assume
that $m(n,1,i)\ge \gamma(n)n.$ We then define
$$
\psi_n^{(1,i)}(f)= \phi_n^{(1,i)}(f)=\oplus_{j=1}^{\gamma(n)}{\rm
diag} ( f(\xi(j,1)), f(\xi(j,2)),...,f(\xi(j,k_j))),
$$
where  $\{\xi(j,1),...,\xi(j,k_j)\}$ is a set of $k_j$ points on the $j$-th
circle of $X_n,$  from $C(X_n)\otimes M_{k(n)}$ into the $i$-th
summand of $B_n.$ We choose $k_j$ so that
$\sum_{j=1}^{\gamma(n)}k_j=m(n,1,i).$ Moreover, we  choose
$\{\xi(j,1),...,\xi(j,k_j)\}$ so that it is $2\pi/n$ dense in (the
$j$-th ) $S^1.$ We may also assume that
$\psi_n^{(1,i)}(1_{C(X_n)\otimes
M_{k(n)}})=\phi_n^{(1,0)}(1_{M_{k(n)}}).$

Now define $\psi_n: A_n\to A_{n+1}$  as follows: Define the
partial map of $\psi_n$ from $C(X_n)\otimes M_{k(n)}$ to
$C(X_{n+1})\otimes M_{k(n+1)}$ to be $\psi_n^{(1,0)}.$ Define the
partial map of $\psi_n$ from $C(X_n)\otimes M_{k(n)}$ to the
$i$-th summand of $B_{n+1}$ to be $\psi_n^{(1,i)}.$ Define
$\psi_n|_{B_n}=\phi_n|_{B_n}. $ Note, here, we identify
$M_{k(n+1)}$ with ${\mathbb C}\otimes M_{k(n+1)}\subset
C(X_{n+1})\otimes M_{k(n+1)}.$

It is standard and easy to verify that $C=\lim_{n\to\infty}(A_n,
\psi_n)$ is a unital simple \CA\, with tracial rank zero (for example,
3.7.8 and 3.7.9 of \cite{Lntext}

Note
that
$$
(K_0(A_n), K_0(A_n)_+, [1_{A_n}])=(K_0(M_{k(n)}\oplus
B_n),K_0(M_{k(n)}\oplus B_n)_+, [1_{M_{k(n)}\oplus B_n}]).
$$
Moreover $(\psi_n)_{*0}=(\phi_n)_{*0}.$ One also have
$$
K_1(A_n)=G^{\gamma(n)}
$$
and $(\psi_n)_{*1}=\kappa_n.$ It follows that
$$
(K_0(C),K_0(C)_+,[1_C], K_1(C))=(K_0(A), K_0(A)_+,[1_A], K_1(A)).
$$
Since $C$ has tracial rank zero and is an AH-algebra (so it
satisfies the UCT), by \cite{Lnmsri}, $C\cong A.$ The lemma follows
from the construction.
\end{proof}

 \begin{thm}\label{MT1}
 Let $A$ be a unital simple A$\T$-algebra with real rank zero and $h: K_1(A)\to
 K_1(A)$ be an isomorphism. Then
 there exists an automorphism $\alpha: A\to A$ which satisfies the
tracial cyclic Rokhlin property such that $\alpha_{*1}=h$ and
$\alpha_{*0}=\id_{K_0(A)}.$
\end{thm}

\begin{proof}
If $K_1(A)=\{0\},$ it suffices to show that there exists an
automorphism on $A$ which has the Rokhlin property. But that
follows from \cite{Phn}. Therefore, for the rest of the proof, we
may assume that $K_1(A)\not=\{0\}.$

Since $K_1(A)$ is tosion free, we may write that
 $K_1(A)=\lim_{n\to\infty}(G_n, \kappa_n),$ where each $G_n$ is
a direct sum of finitely many copies of $\Z$ and $\kappa_n$ are
injective.

Let $h: K_1(A)\to K_1(A)$ be the given isomorphism and let
$h^{-1}$ be the inverse.  We may assume, by passing to a subsequence if necessary,
 that there  are
homomorphisms $h_n, {\bar h}_n: G_n\to G_{n+1}$ such that
$$
h\circ \kappa_{n,\infty}=\kappa_{n+1,\infty}\circ h_n
\andeqn
h^{-1}\circ \kappa_{n,\infty}=\kappa_{n+1,\infty}\circ {\bar
h}_n
$$
on $G_n.$
 By passing to a subsequence if necessary, we may assume
that,
\begin{eqnarray}\label{adTp1}
h_{n+1}\circ \kappa_n=\kappa_{n+1}\circ h_n, \,\,\, {\bar
h}_{n+1}\circ h_n=\kappa_{n+1}\circ \kappa_n \andeqn h_{n+1}\circ
{\bar h}_n =\kappa_{n+1}\circ \kappa_n.
\end{eqnarray}

Suppose that $G_n$ is $c(n)$ copies of $\Z.$ Homomorphisms
 $h_n,$ $h_n^{-1}$ and $\kappa_n$ can be represented by
matrices  $(a_{i,j}^{(n)}),$ $(b_{i,j}^{(n)})$ and
$(c_{i,j}^{(n)})$ where $a_{i,j}^{(n)},$ $b_{i,j}^{(n)}$ and
$c_{i,j}^{(n)}$ are integers.

Let
$J(n)=\max\{|a_{i,j}^{(n)}|+|b_{i',j'}^{(n)}|+|c_{i'',j''}^{(n)}|:
\,i,j,i',j',i'',j''\}.$

We identify $X_n$ with $\gamma_n$ copies of the unit circle with a
common point at 1. Denote $u(n,i)\in C(X_n)$ the function on $X_n$
which is the identity map on the $i$th circle and 1 everywhere
else.

 Define
$$
\beta_n(u(n,i))=\prod_{j=1}^{\gamma(n+1)}z(n+1,j,i),
$$
where $z(n+1,j,i)=u(n+1,j)^{a_{i,j}^{(n)}}, $ if
$a_{i,j}^{(n)}\not=0,$ and $z(n+1,j,i)=u(n,i)(1)$ if
$a_{i,j}^{(n)}=0,$ $i=1,2,...,\gamma(n).$

Note $\beta_n$ gives a homomorphism from $C(X_n)\to C(X_{n+1})$
such that $(\beta_n)_{*1}=h_n.$

Define
$$
{\bar \beta_n}(u(n,i))=\prod_{j=1}^{\gamma(n+1)}z'(n+1,j,i),
$$
where $z'(n+1,j,i)=u(n+1,j)^{b_{i,j}^{(n)}}, $ if
$b_{i,j}^{(n)}\not=0,$ and $z'(n+1,j,i)=u(n,i)(1)$ if
$b_{i,j}^{(n)}=0,$ $i=1,2,...,\gamma(n).$


 We now write
$A=\lim_{n\to\infty}(A_n, \phi_n),$ where $A_n=C(X_n)\otimes
M_{k(n)}\oplus B_n$ which satisfying (1)-(6) in 
Lemma \ref{smallpro2} with $a_n=J(n)n$ We
will also use $\beta_n: A_n\to A_{n+1}$ for the \hm\, defined by
$(\beta_n)|_{C(X_n)\otimes M_{k(n)}}=\beta_n\otimes {\rm
id}_{M_{k(n)}}$ and $(\beta_n)|_{B_n}=\phi_n|_{B_n}.$

Let ${\cal F}_n\subset {\cal F}_{n+1}$ be a sequence of finite
subsets of $A$ such that the union of those finite subsets is
dense in $A.$ Without loss of generality, we may assume that
${\cal F}_n\subset \phi_{n, \infty}({\cal G}_n),$ where ${\cal
G}_n$ is a finite subset of $A_n.$

It follows from \cite{Phn} that there is an approximately inner
automorphism $\sigma\in Aut(A)$ which has the tracial Rokhlin
property. It follows from \cite{LO} that $\sigma$ has tracial
cyclic Rokhlin property.

Since $\sigma$ is approximately inner, by passing to a subsequence if necessary, by Lemma \ref{smallpro2},
we may assume that, there is a unitary $v_n\in A_{n+1}$ such that
\begin{eqnarray}\label{appinn1}
\sigma\circ \phi_{n,\infty}(f)\approx_{1/2^{n+2}} {\rm ad}\,\phi_{n+1,\infty}(v_n)\, \circ \phi_{n, \infty}(f)
\andeqn
\end{eqnarray}
\begin{eqnarray}\label{appinn2}
\sigma^{-1}\phi_{n,\infty}(f)\approx_{1/2^{n+2}}{\rm ad}\, \phi_{n+1, \infty}(v_n^*)\circ
\phi_{n, \infty}(f)
\,\,\,{\rm on}\,\,\,
{\cal G}_n.
\end{eqnarray}
Define  $\sigma_n: A_{n+1}\to A_{n+1}$ by $\sigma_n(f)={\rm ad}\, v_n(f)$
and $\dt_n: A_{n+1}\to A_{n+1}$ by $\dt_n(f)={\rm ad}\, v_n^*(f)$ for $f\in A_n.$

Since $a_n=J(n)n,$ by passing to a subsequence if necessary, we
may assume that $\phi_n^{(1,1)}: C(X_n)\otimes M_{k(n)}\to
C(X_{n+1})\otimes M_{k(n+1)}$ has the following form:
$$
\phi_n^{(1,1)}(f)={\rm diag}(f\circ s_n,
f(\zeta(n,1)),f(\zeta(n,2),...,f(\zeta(n,t(n))),
$$
where $s_n$ is a standard map (see Definition \ref{Standard}) and
$\{\zeta(n,1)),\zeta(n,2),...,\zeta(n,t(n)\}$ contains subsets
$\{\xi_1,...,\xi_N\},$ $\{\xi_{N+1},....\xi_{2N}\}$,...,$\{\xi_{(\gamma(n)-1)N+1},...,
\xi_{\gamma(n)N}\}$ such that
$\{\xi_{(j-1)N+1},...,\xi_{jN}\}$ evenly divide the $j$-th copy of $S^1$ and
$N=\Lambda(n, {\cal G}_n).$

 Define a monomorphism $\omega_n: C(X_n)\otimes
M_{k(n)}\to A_{n+2}$ as follows:
%
Define
$$
\omega_n^{(1,1)}(f)=
\beta_{n+1}(g_1)+\phi_{n+1}(g_2),
$$
where $g_1=\sigma_{n}({\rm diag}(f\circ s_n, 0\cdots 0))$ and
$g_2=\sigma_{n}(0, f(\zeta(n,1)),f(\zeta(n,2)),...,f(\zeta(n,t(n))).$
Define $\omega_{n}^{(1,i)}=\phi_{n+1}^{(1,i)}\circ \phi_{n}|_{C(X_n)\otimes M_{k(n)}}.$ Let $E_{n+2}=
\phi_{n,n+2}(d_n),$ where $d_n$ is the identity of $C(X_n)\otimes
M_{k(n)}.$ Denote by
$$
E_{n+1}'={\rm diag}(1,\overbrace{0,...,0}^{t(n)})\in C(X_{n+1})\otimes M_{k(n+1)}.
$$

Define a monomorphism ${\bar\omega}_n: C(X_n)\otimes M_{k(n)}\to
E_{n+2}A_{n+2}E_{n+2}$ as follows:

Define
$$
{\bar \omega}_n^{(1,1)}(f)=
{\bar \beta_{n+1}}(g_1')+\phi_{n+1}(g_2'),
$$
where $g_1'=\dt_{n}({\rm diag}(f\circ s_n, 0\cdots 0))$ and
$g_2'=\dt_{n}(0, f(\zeta(n,1)),f(\zeta(n,2),...,f(\zeta(n,t(n))).$
Define ${\bar \omega}_n^{(1,i)}=\phi_{n+1}^{(1,i)}\circ \phi_n|_{C(X_n)\otimes M_{k(n)}}.$

For any unitary $u\in A_{4n+3},$ define $\Phi: C(X_{4n+1})\otimes M_{k(4n+1)}\to A_{4n+5}$ by
$$
\Phi(f)={\bar \bt}_{4n+4}\circ \dt_{4n+3}\circ\phi_{4n+3}\circ {\rm ad}\, u\circ \bt_{4n+2}(g_1).
$$
Since
$$
({\bar\bt}_{4n+4}\circ \dt_{4n+3} \circ\phi_{4n+3}\circ {\rm ad}\, u\circ \bt_{4n+2}\circ \phi_{4n+1})_{*1}=(\phi_{4n+1, 4n+5})_{*1},
$$
$$
(\Phi)_{*1}=(\phi_{4n+2,4n+5}\circ E'_{4n+1}\phi_{4n+1}E'_{4n+1})_{*1}.
$$

By choice of  $\Lambda(n, {\cal G}_n),$  by (\ref{adTp1}) and by
applying Lemma \ref{LLambda}, we obtain unitaries $u_k\in
E_{k}A_{k}E_{k}$ such that (with $u_1=d_1$)
\begin{eqnarray}\label{inter1}
{\rm ad}\, u_{4(n+1)+1}\circ {\bar \omega}_{4n+3}\circ {\rm ad}\,
u_{4n+3}\circ \omega_{4n+1}(f)\approx_{1/2^n} \phi_{4n+1,4(n+1)+1}(f)
\,\,\, {\rm on}\,\,\, d_{4n+1}{\cal G}_{4n+1}
\end{eqnarray}
as well as
\begin{eqnarray}\label{inter2}
{\rm ad}\, u_{4(n+1)+3}\circ \omega_{4(n+1)+3}\circ {\rm ad}\,
u_{4n+3}\circ {\bar \omega}_{4n+3} (f)\approx_{1/2^{n}} \phi_{4n+3,
4(n+1)+3}(f) \,\,\, {\rm on}\,\,\, d_{4n+3}{\cal G}_{4n+3}
\end{eqnarray}
($n=0,1,2,...$)

 Define $\omega_{4n+1}'={\rm ad}\,u_{4n+3}\circ \omega_{4n+1}$
and ${\overline{\omega}'}_{4n+3}={\rm ad}\, u_{4n+5}\circ
{\bar\omega}_{4n+3}.$

Now define $\alpha_n: A_{4n+1}\to A_{4n+3}$ by
$$
(\alpha_n)|_{C(X_{4n+1})\otimes M_{k(4n+1)}}=\omega_{4n+1}'\andeqn
(\alpha_n)|_{B_{4n+1}}=(\phi_{4n+2}\circ \sigma_{4n+1}\circ
\phi_{4n+1})|_{B_{4n+1}}
$$
and define ${\bar\alpha}_n: A_{4n+3}\to A_{4n+5}$ by

$$
({\bar \alpha}_n)|_{C(X_{4n+3}) \otimes
M_{k(4n+3)}}={\overline{\omega'}_{4n+3}\andeqn ({\bar
\alpha}_n)|_{B_{4n+3}}=(\phi_{4n+4}\circ \dt_{4n+3}\circ
\phi_{4n+3})|_{B_{4n+3}}}
$$

From (\ref{appinn1}) and (\ref{appinn2}), we have that
\begin{eqnarray}\label{inter3}
{\bar \af}_{n+1}\circ \af_n(f)\approx_{1/2^{n+1}} \phi_{4n+1, 4(n+1)+1}(f)\,\,\,{\rm on}\,\,\,f\in (1-d_{4n+1}){\cal G}_{4n+1}
\end{eqnarray}

\begin{eqnarray}\label{inter4}
{\af}_{n+1}\circ {\bar \af}_n(f)\approx_{1/2^{n+1}} \phi_{4n+3, 4(n+1)+3}(f)\,\,\,{\rm on}\,\,\,f\in (1-d_{4n+1}){\cal G}_{4n+1}
\end{eqnarray}

It follows from
(\ref{inter1}), (\ref{inter2}), (\ref{inter3}) and
(\ref{inter4}), one has the following approximately
intertwining:
$$
\begin{array}{ccccccccccc}
A_1& {\stackrel{\phi_{1,5}}\longrightarrow}  &A_5  &{\stackrel{\phi_{5,9}}\longrightarrow} & A_9 & {\stackrel{\phi_{9,13}}\longrightarrow}&A_{17}\cdots\\
 &\searrow_{\af_1}\,\, \nearrow_{{\bar \af}_1} & & \searrow_{\af_2}\,\, \nearrow_{{\bar \af}_2} &&\searrow_{\af_3} \,\,\,\nearrow_{{\bar \af}_3}\\
&  A_3 &{\stackrel{\phi_{3,7}}\longrightarrow}& A_7 &{\stackrel{\phi_{7, 11}}\longrightarrow}&A_{11} &
{\stackrel{\phi_{11, 15}}\longrightarrow}\cdots
\end{array}
$$

It follows from the Elliott approximately intertwining argument that
$\{\alpha_n\}$ and $\{{\bar\alpha}_n\}$  define two automorphisms
$\alpha$ and $\alpha^{-1}.$

Moreover, one checks that
$$
\alpha_{*0}={\rm id}_{K_0(A)}\andeqn \alpha_{*1} =h.
$$

It remains to show that $\alpha$ has tracial cyclic Rokhlin
property.

Fix $\ep>0,$ $m>0$ and a finite subset ${\cal F}\subset A.$  Since $\sigma$
has tracial cyclic Rokhlin property, there are mutually orthogonal
projections $e_1, e_2,...,e_m$ such that (with $e_{m+1}=e_1$)
\begin{eqnarray}\label{adTp2}
\|e_ia-ae_i\|<\ep/8,\,\,\,i=1,2,...,m, a \in {\cal F}
\end{eqnarray}

\begin{eqnarray}\label{adTp3}
\|\sigma(e_i)-e_{i+1}\|<\ep/8,\,\,\,i=1,2,...,m
\end{eqnarray}
and
\begin{eqnarray}\label{adTp4}
\tau(1-\sum_{i=1}^me_i)<\ep/8
\end{eqnarray}

for all $\tau\in T(A).$

Since $\sigma$ is approximately inner, there is a unitary $v\in A$ such that
$$
\|v^*e_iv-e_{i+1}\|<\ep/8,\,\,\,i=1,2,...,m.
$$

Without
loss of generality, we may assume that there is a finite subset
${\cal G}\subset A_{4n+1}$ such that ${\cal F}\subset \phi_{4n+1,
\infty}({\cal G})$  for some sufficiently large $n.$

It follows from
the standard argument 
(See \cite[Proposition L.2.2]{Wo} for example.) that for any $\dt>0,$
(for sufficiently large
$n$) there are
 projections $q_i\in A_{4n+1}$
 such that
 \begin{eqnarray}\label{adTP9}
 \|e_i-\phi_{4n+1,\infty}(q_i)\|<\min\{\dt, \ep/8\}
 ,\,\,\,i=1,2,...,m.
\end{eqnarray}
Since $\phi_{n,\infty}$ is assumed to be injective, with
sufficiently small $\dt,$ 
from \cite[Lemma 2.5.6]{Lntext} we may assume that
$q_1,q_2,...,q_m$ are mutually orthogonal. Let ${\bar
q}_i=\phi_{4n+1,\infty}(q_i).$ We may also assume that there
is a unitary $v'\in A_{4n+1}$ such that
$\phi_{4n+1,\infty}(v')=v.$

With even larger $n,$ we may assume that, without loss of
generality,
\begin{eqnarray}\label{adTp10}
\alpha^{-i}({\cal F}),\, \alpha^{-i}({\bar q_i})\in
\phi_{4n+1,\infty}(A_{4n+1}), \,\,\, i=1,2,...,m.
\end{eqnarray}

Note that $(1-d_{4n+1})$ commutes with every element in
$A_{4n+1}.$
Put
$p_i={\bar q}_i\alpha^{i-1}(\phi_{4n+1,\infty}(1-d_{4n+1})),$ $i=1,2,...,m.$
Since, by (\ref{adTp10}),
$$
\alpha^{i-1}\circ \phi_{4n+1,\infty}(1-d_{4n+1}) {\bar q_i}=
\alpha^{i-1}(\phi_{4n+1,\infty}(1-d_{4n+1}) \alpha^{-i+1}({\bar
q_i}))
$$
$$
=\alpha^{i-1}(\alpha^{-i+1}({\bar
q_i})\phi_{4n+1,\infty}(1-d_{4n+1}))= {\bar q_i}\alpha^{i-1}\circ
\phi_{4n+1,\infty}(1-d_{4n+1}),
$$
$\{p_1,p_2,...,p_m\}$ are mutually orthogonal projections.

By (\ref{adTp10}), if $a\in {\cal F},$
$$
p_ia={\bar
q}_i\alpha^{i-1}(\phi_{4n+1,\infty}(1-d_{4n+1})\alpha^{-i+1}(a))
={\bar q}_ia\alpha^{i-1}(\phi_{4n+1,\infty}(1-d_{4n+1})).
$$
Therefore, by (\ref{adTP9}) and (\ref{adTp2}), we have

(1) $\|p_ia-ap_i\|
=\|{\bar q}_ia\alpha^{i-1}((\phi_{4n+1,\infty}(1-d_{4n+1})))-ap_i\|
<\ep/8+\ep/8+\ep/8<\ep,\,\,\,i=1,2,...,m.$

Since
$$
\alpha(p_i)=
\alpha({\bar q}_i\alpha^{i-1}(\phi_{4n+1,\infty}(1-d_{4n+1})))
= \alpha({\bar q}_i)\alpha^{i}(\phi_{4n+1,\infty}(1-d_{4n+1}))),
$$
by (\ref{adTp2}) and (\ref{adTp3}),

(2) $\|\alpha(p_i)-p_{i+1}\|=\|(\alpha({\bar q}_i)-{\bar q_{i+1}})
\alpha^i(\phi_{4n+1,\infty}(1-d_{4n+1}))\|<\ep,$ $i=1,2,...,m.$

Note, since
$$
\lim_{n\to\infty}\sup\{\tau(\phi_{n,\infty}(d_n)):\tau\in
T(A)\}=0,
$$
we may assume that
$$
\tau(\phi_{4n+1,\infty}(d_{4n+1}))<\ep/4m
$$
for all $\tau\in T(A).$

Therefore, using the fact that

 $\alpha_{*0}={\rm id}_{K_0(A)},$
$$
\tau(p_i)=\tau({\bar
q}_i\af^{i-1}(\phi_{4n+1,\infty}(1-d_{4n+1})))> \tau({\bar
q}_i)(1-\ep/4m),\,\,\,i=1,2,...,m
$$
Note also, by (\ref{adTP9})  and with $\dt<1,$
$$
\tau({\bar q}_i)=\tau(e_i),\,\,\,i=1,2,...,m.
$$

Applying (\ref{adTp4}), we also have

(3)
\begin{align*}
\tau(1-\sum_{i=1}^m p_i) &< 1- (1 - \frac{\ep}{4m})\sum_{i=1}^m\tau({\bar q}_i)\\
&= 1- (1 - \frac{\ep}{4m})\sum_{i=1}^m\tau(e_i)\\
&= (1 - \sum_{i=1}^m\tau(e_i)) + \frac{\ep}{4m}\sum_{i=1}^m\tau(e_i)\\
&< \frac{\ep}{8} + \frac{\ep}{4}\tau(e_1)
< \frac{\ep}{8} + \frac{\ep}{4} < \ep
\end{align*}
 for all $\tau\in T(A).$

\end{proof}

\begin{Cor}\label{MT1C}
Let $A$ be a unital simple $A\T$-algebra with real rank zero and
let $h: K_1(A)\to K_1(A)$ be an isomorphism. Suppose that $\af$ is
the automorphism given by Theorem \ref{MT1}. Then the crossed
product $A\rtimes_{\af}\Z$ has tracial rank zero and can be
written as an AH-algebra.
\end{Cor}

\begin{proof} Since $\af_{*0}={\rm id}_{K_0(A)}$ and has the
tracial cyclic Rokhlin property, by  \cite{Lngpots}
$A\rtimes_{\af}\Z$ has tracial rank zero. It follows from the
classification theorem \cite{Lnmsri} that $A\rtimes_{\af}\Z$ is in
fact simple AH-algebra.
\end{proof}

\section{Minimal dynamical systems and non-commutative dynamical systems}

It seems that a more general construction based on 
Theorem \ref{MT1} can be made
so that Theorem \ref{MT1} can covers  all simple AH-algebras with real rank zero and stable rank one.
However, in this section, we present a different type of construction. If one believes that automorphisms
on \CA s with Rokhlin
property  is closely related to minimal homeomorphisms on compact metric spaces, then the construction
that we presented here sheds some light on these relation.
Given a connected finite CW complex $X,$  there is a standard way to construct a unital simple \CA\, with similar
topological information:

\vspace{0.2in}

\begin{Prop}\label{NCX}
Let $X$ be a connected finite CW-complex space and let $D$ be
a dimension group which is a countable dense subgroup of $\Q.$
Let $\rho_X: K_0(C(X))\to \Z$ be the dimension map.
Then there exists a unital separable $AH$-algebra $A_X$ with tracial rank zero and
with a unique tracial state $\tau$ such that
$$
(K_0(A_X), K_0(A_X)_+)=D\oplus {\rm ker}{ \rho}_X,\,\, {\rm ker}\rho_{A_X}={\rm ker}\rho_X\andeqn
K_1(A_X)=K_1(C(X)).
$$
Moreover, there is a unital embedding $j: C(X)\to A_X$ such that
$$
j_{*0}|_{{\rm ker} \rho_X}={\rm id}_{{\rm ker} \rho_X}\andeqn
j_{*1}={\rm id}_{K_1(C(X))}.
$$
\end{Prop}
(See \cite{Lnact} for the construction).

Suppose that there is a minimal homeomorphism $\psi:X\to X.$ One can then construct
an automorphsim $\af$ on $A_X$ associated with $\psi.$

\begin{thm}\label{XR}
Let $X$ be a connected finite CW-complex space,
and let $\psi\colon X \rightarrow X$ be  a minimal homeomorphism.
Denote by $\psi^{\natural}: C(X)\to C(X)$ the automorphism defined by $\psi^{\natural}(f)=f\circ \psi.$
Then there is an automorphism $\af: A_X\to A_X$ which satisfies the tracial cyclic
Rokhlin property such that
$$
\af_{*0}|_{{\rm ker} \rho_X}=(\psi^{\natural})_{*0}|_{{\rm ker} \rho_X}\andeqn
\,\,\, \af_{*1}=(\psi^{\natural})_{*1}.
$$
\end{thm}

Theorem \ref{XR} is a corollary of a more general Theorem \ref{N2MT} below.

\begin{Lem}\label{Marr1}
Let $X$ be a compact connected metric space and let $\psi: X\to X$
be a minimal homeomorphism with a $\psi$-invariant Borel
probability measure $\mu.$ Let $S$ be an infinite subset of
positive integers. Then, for any $\ep>0,$ there is an integer
$n\in S$ and there is a finite $\ep$-dense set $\{x_1, x_2,
...,x_n\}\subset X$ such that for any closed subset $F\subset X,$
\begin{eqnarray}\label{Marr}
\mu_1(F)\le \mu_2(F_{\ep})\andeqn \mu_2(F)\le \mu_1(F_{\ep}),
\end{eqnarray}
where $\mu_1$ is the measure concentrated on $\{x_1,x_2,...,x_n\}$
with $\mu_1(\{x_j\})={1\over{n}}$ and $\mu_2$ is the measure
concentrated on $\{\psi(x_1), \psi(x_2),...,\psi(x_n)\}$ with
$\mu_2(\{\psi(x_j)\})={1\over{n}}.$
\end{Lem}

\begin{proof} First we note that any open neighborhood $O(x)$
of any $x\in X$ must have positive $\mu$ measure. Otherwise,
assume that $\mu(O(x))=0.$ Put
$$
O=\cup_{n\in \Z} \psi^n(O(x)).
$$
Then $O$ is open and $\psi(O)=O.$ It follows that $F=X\setminus O$
is invariant under $\psi.$ The minimality of $\psi$ implies that
$F=\emptyset.$ This implies that $O=X.$ But $\mu(O)=0.$ A
contradiction.

Fix $\dt>0.$ Since $X$ is compact, by what we have just proved, it
is easy to find finitely many disjoint Borel subsets
$G_1,G_2,...,G_m$ such that $\mu(G_i)>0$ and ${\rm diam}(G_i)<\dt$
for $i=1,2,...,m,$ and
$$
\cup_{i=1}^m G_i=X.
$$
Let $a=\min\{\mu(G_i): 1\le i\le m\}.$ Choose $n\in S$ so that
$$
{1\over{n}}<{a\over{2m^{m+1}}}.
$$
There is an integer $0<k_i<n$ such that
$$
{k_i\over{n}}\le \mu(G_i)<{k_i+1\over{n}},i=1,2,...,m.
$$
By minimality of $\psi,$ no single point has positive
$\mu$-measure. Thus each $G_i$ contains infinitely many points.
Choose $k_i$ points in $G_i,$ $i=1,2,...,m-1.$
Then
$$
\sum_{i=1}^m k_i\ge n-m.
$$
Choose $k_m'=k_m+(n-\sum_{i=1}^m k_i)$ points in $G_m.$
Note that
\begin{eqnarray}\label{Marr1+}
{k_m'\over{n}}-\mu(G_m)\le {m\over{n}}.
\end{eqnarray}

 We obtain
$n$ points $x_1,x_2,...,x_n$ in $X.$ Now define $\mu_{\dt,1}$  to be the
measure concentrated on $\{x_1,x_2,...,x_n\}$ with
$\mu_{\dt,1}(\{x_j\})={1\over{n}}$ and define $\mu_{\dt,2}$ to be
the measure concentrated on $\{\psi(x_1), \psi(x_2),...,\psi(x_n)\}$
with $\mu_{\dt,2}(\{\psi({x_j})\}={1\over{n}}.$ It is clear that when $\dt\to
0,$ $m\to\infty.$ Fix $\xi(\ep,j)\in G_j.$ Note that for  each
$f\in C(X),$
$$
\sum_{j=1}^mf(\xi(\dt,j))\mu(G_j)\to \int_Xfd \mu
$$
as $\dt\to 0.$ It follows  that
$$
\int_X f d\mu_{\dt,1}\to \int_X f d\mu
$$
for all $f\in C(X)$ as $\dt\to 0.$ We also have
$$
\int_X f\mu_{\dt,2}=\int_X f\circ \psi d\mu_{\dt,1}\to \int_X
f\circ \psi d\mu=\int_X f d\mu
$$
for all $f\in C(X).$

Now let $\ep>0$ and $y_1,y_2,...,y_K$ be a subset of $X$ such that
$$
\cup_{j=1}^K O(y_j,\ep/4)=X.
$$

Let ${\cal O}$ be the finite collection of all possible union of
$O(y_j,\ep/4).$ Let ${\cal O}'$ be the finite collection of all
possible union of $O(y_j,\ep/2).$ If $G\subset {\cal O}$ we denote
by $G'$ the corresponding union of $O(y_j,\ep/2)$'s.

Let
$$
b=\min\{\mu(G')-\mu(G): G\subset {\cal O} \andeqn G\not=X\}.
$$
The connectedness of $X$ together the first paragraph of this
proof imply that $b>0.$

Choose $\dt$ sufficiently small so that, for all $G\subset {\cal
O}$
$$
\mu_{\dt,i}(G)-{b\over{4}}\le \mu(G))\le
\mu_{\dt,i}(G)+{b\over{4}}
$$
and, for any $G'\in {\cal O}',$
$$
\mu_{\dt,i}(G')-{b\over{4}}\le \mu(G')\le
\mu_{\dt,i}(G')+{b\over{4}}
$$
$i=1,2.$ Thus, if $G\not=X,$
$$
\mu_{\dt,i}(G)+{b\over{4}}\le \mu(G)+{b\over{4}}\le
\mu(G')-{b\over{4}}\le \mu_{\dt,i}(G'),i=1,2.
$$

Now let $F$ be a closed subset of $X.$ If $F_{\ep}=X,$ then
(\ref{Marr}) follows. Now suppose that $F_{\ep}\not=X.$

 Let
$$
G=\cup_{F\cap O(y_j,\ep/4)\not=\emptyset} O(y_j,\ep/4) \subset
{\cal O}
$$
 Then $G\not=X.$ It follows that
$$
\mu_{\dt,1}(F)\le \mu(G)\le \mu_{\dt,2}(G)+{b\over{4}}\le
\mu_{\dt,2}(G')\le \mu_{\dt,2}(F_{\ep}).
$$
Similarly
$$
\mu_{\dt,2}(F)\le \mu_{\dt,1}(F_{\ep}).
$$
\end{proof}

\begin{Lem}\label{EVEN}
Let $X$ be a compact connected metric space (with infinitely many
points) and let $\psi: X\to X$ be a minimal homeomorphism. Let $S$
be an infinite subset of $\N.$ Then, for any $\ep>0,$ there is
$n\in S$ and there is a finite $\ep$-dense set $\{x_1, x_2,
...,x_n\}\subset X$ such that there is a permutation $s:
(1,2,...,n)\to (1,2,...,n)$ such that
$$
{\rm dist}(x_j, \psi(x_{s(j)}))<\ep, j=1,2,...,n.
$$
\end{Lem}

\begin{proof}
Let $\mu$ be a $\psi$-invariant Borel probability measure.  Let
$\{x_1,x_2,...,x_n\}$ be a finite subset satisfying the conclusion
of Lemma \ref{Marr1}.

For any subset $S_1$ of $l$ elements in $\{x_1,x_2,...,x_n\},$
choose a closed subset $F\cap \{x_1,x_2,...,x_n\}=S_1$ and
$$
F\subset\{x\in X: {\rm dist}(x,\{x_1,x_2,...,x_n\})<\ep/4\}.
$$
Then
$$
\mu_1(F)\le \mu_2(F_{\ep/4}).
$$
Thus $F_{\ep/4}$ contains at least $l$ elements in
$\{\psi(x_1),\psi(x_2),...,\psi(x_n)\}.$ Therefore there is a
subset $S_2$ of $l$ elements in
$\{\psi(x_1),\psi(x_2),...,\psi(x_n)\}$ such that, for any element
$\xi\in S_1$ there is $\xi'\in S_2$ such that
$$
{\rm dist}(\xi, \xi')<\ep.
$$
By the "Marriage Lemma" (see \cite{HV}), there is a permutation $s$
which meets the requirement.

\end{proof}


\begin{Lem}\label{Lem1.2}
Let $X$ be a compact connected metric space (with infinitely 
many points) with a metric $d$,
and $\psi\colon X \rightarrow X$ be a minimal homeomorphism.
For any $\epsilon > 0$ and a finite set ${\mathcal F} \subset C(X)$
there exists $\delta > 0$ such that
if $\{x_i\}_{i=1}^n$ are points in $X$ which satisfy the property that 
there is a permutation $s\colon (1,2,\dots,n)\rightarrow (1,2,\dots,n)$ 
such that 
$$
{\rm dist}(x_j,\psi(x_{s(j)})) < \delta,
$$
then 
there is a permutation $n \times n$ matirx $U$
such that
$$
\|{\rm diag}(f(x_1),\dots, f(x_n)) -
U{\rm diag}(f(\psi(x_1)), \dots, f(\psi(x_n))U^*\| < \epsilon,\ \
\forall f \in {\mathcal F}.
$$
\end{Lem}

\begin{proof}
Let $\ep > 0$ and ${\mathcal F} \subset C(X)$ be a finite set.
Then there exists $\dt > 0$ such that 
if $d(x,x') < \dt$ then $|f(x) - f(x')| < \ep$ for any $f \in {\mathcal F}$.

Let $\{x_1, x_2, \dots, x_n\}$ be points in $X$ 
satisfy that 
there is a permutation $s\colon (1,2,\dots,n)\rightarrow (1,2,\dots,n)$ 
such that 
$$
{\rm dist}(x_j,\psi(x_{s(j)})) < \delta.
$$

Take a permutation matrix $U$ correspondent to $s$, then 
\begin{align*}
&\|{\rm diag}(f(x_1),\dots, f(x_n)) -
U{\rm diag}(f(\psi(x_1)), \dots, f(\psi(x_n))U^*\|\\
&= 
\|{\rm diag}(f(x_1),\dots, f(x_n)) -
{\rm diag}(f(\psi(x_{s(1)})), \dots, f(\psi(x_{s(n)}))\|\\
&< \ep
\end{align*}
\end{proof}

\begin{Lem}\label{Nlem1.2}
Let $X$ be a compact connected metric space with a minimal
homeomorphism $\psi: X\to X.$ Let $\{f_n\}$ be a dense sequence in
$C(X),$ let $\{\ep_n\}$ be a sequence of positive numbers and let
$\{k(n)\}$ be a sequence of positive integers. Then there exists a
sequence of positive numbers $\{\dt_n\}$ which has the following
property: If $\{x(1,n),..., x(l(n),n)\}$ are points in $X$ which
satisfy the following property: there is a permutation $s_n:
(1,2,...,l(n)) \to (1,2,...,l(n))$ such that
$$
{\rm dist}(x(j,n),\psi(x(s_n(j),n)))<\dt_n
$$
then there is a sequence of permutation $l(n)\times l(n)$ matrices
$U_n$ such that
$$
W_n^*{\rm diag}(f(x(1,n)),f(x(2,n)),...,f(x(l(n),n)))W_n
$$
$$
\approx_{\ep_n} {\rm
diag}(f(\psi(x(1,n))),f(\psi(x(2,n))),...,f(\psi(x(l(n),n))))
$$
for all $f\in \{(a_{i,j})\in C(X)\otimes M_{k(n)} : a_{i,j}\in
\{f_1,f_2,...,f_n\}\cup{\mathbb C}\},$ where $W_n=1_{k(n)}\otimes U_n$ is a matrix in
$M_{k(n)l(n)}.$
\end{Lem}

\begin{proof}
This follows from  Lemma \ref{Lem1.2} immdiately.
\end{proof}

We shall construct a model automorphism with tracial Rokhlin property
on some AH algebras which contains the class of $A{\mathbb T}^n$-algebras
and $AI$-algebras.

Let $X$ be a compact metric space with infinitely many points. Let
$(X,\psi)$ be a minimal dynamical system. If $(X,\psi)$ is
uniquely ergodic, then $(X,\psi)$ has mean dimension zero. If $X$
has finite covering dimension, then $(X,\psi)$ always has mean
dimension zero (\cite{LW}).

\begin{Lem}\label{Lrok}
Let $X$ be a compact metric space and let $\psi: X\to X$ be a
minimal homeomorphism. Suppose that $(X,\psi)$ has mean dimension
zero. Then, for any integer $N\ge 1,$ any $\eta>0$ and  any
$\dt>0,$ there are mutually disjoint open subsets
$G_1,G_2,...,G_L$ such that

(1) ${\rm diam}(G_i)<\eta,$ $i=1,2,...,L,$

(2) $\psi^i(\cup_{s=1}^LG_s)\cap
\psi^j(\cup_{s=1}^LG_s))=\emptyset,$ if $i\not=j$ and
$i,j=0,1,2,...,N-1,$

(3) $\mu(\cup_{j=0}^{N-1}\psi^j(\cup_{i=1}^LG_i))>1-\dt$ for all $\mu\in
T_{\psi}$ and

(4) $\mu(\partial(G_i))=0$ for $i = 1, 2, \dots, L$ and
all $\mu\in T_{\psi},$ where
$T_{\psi}$ is the set of all $\psi$-invariant Borel probability
measures
\end{Lem}

This basically follows from Lemma 3.4 and Lemma 3.5 of \cite{Lrokm}.

\begin{Lem}\label{LDig}
Let $X$ be a compact metric space and let $\ep>0.$ Suppose that
${\cal G}\subset C(X)$ satisfies the following: if ${\rm
dist}(x,x')<\dt,$ then
$$
|f(x)-f(x')|<\ep/m^2
$$
for all $f\in {\cal G}.$ Denote by ${\cal F}=\{(a_{i,j})\in
C(X)\otimes M_m: a_{i,j}\in {\cal G}\}.$ Let $\phi: C(X)\otimes M_m\to C(X)\otimes
M_{m(k+l)}$ be defined by
$$
\phi(f)={\rm
diag}(\overbrace{f,f,...,f}^l,f(\xi_1),f(\xi_2),...,f(\xi_k))
$$
for all $f\in C(X)\otimes M_m.$ Then if ${\rm dist}(x,x')<\dt,$
$$
\|\phi(f)(x)-\phi(f)(x')\|<\ep
$$
for all $f\in {\cal F},$ where the norm is operator norm on
$M_{m(k+l)}.$
\end{Lem}

\begin{Dfn}\label{Construction}
{\rm Let $X$ be a compact connected metric space with infinitely
many points. Let $d_X: K_0(C(X))\to \Z$ be the dimension map.
Write $K_0(C(X))=\Z\oplus G_{00}$ and $K_1(C(X))=G_1,$ where
$G_{00}={\rm ker} d_X.$ Fix two sequences $\{a_n\}$ and $\{b_n\}$ and a sequence
$\{k_n\}$ of positive integers (with $k_n\to \infty$ as
$n\to\infty$). We assume that $k_1=1$ and $k_{n+1}=b_nk_n,$ where
$b_n$ a sequence of positive numbers such that $b_n\ge 2.$ Let
$H_{00}=\lim_{n\to\infty}(G_{00}, \kappa_{n,n+1}^{(00)}),$ where
$\kappa_{n,n+1}^{(00)}: G_{00}\to G_{00}$ is defined by
$\kappa_{n,n+1}^{(00)}(g)=a_n g$ for all $g\in G_{00}.$ Let $D$ be
the dimension group defined by $D=\lim_{n\to\infty} (\Z,
\kappa_{n,n+1}^{(d)}),$ where $\kappa_{n,n+1}^{(d)}(g)=b_n g$ for all
$g\in \Z.$ 
Define $H_0=D\oplus H_{00}=\lim_{n,\infty}(\Z\oplus G_{00}, \kappa_{n,n+1}^{(0)}),$
where $\kappa_{n,n+1}^{(0)}=\kappa_{n,n+1}^{(d)}\oplus \kappa_{n,n+1}^{(00)}.$ 
Define $H_1=\lim_{n\to\infty} (G_1,
\kappa_{n,n+1}^{(0)}),$ where $\kappa_{n,n+1}^{(1)}(g)=a_n g$ for
all $g\in G_1.$ By passing to a subsequence, without changing
$H_{00},$ $D$ or $H_,$ one may assume that $b_n=a_n+l(n)$ and
$$
\lim_{n\to\infty}{a_n\over{b_n}}=0.
$$
There is a unital simple AH-algebra $A=A(X,\{a_n\}, \{b_n\})$ such
that
$$
(K_0(A), K_0(A)_0,[1_A], K_1(A))= (D\oplus H_{00}, (D\oplus H_{00})_+,1,
H_1),
$$
where
$$
(D\oplus H_{00})_+=\{(d,x): d>0, d\in D, x\in
H_{00}\}\cup\{(0,0)\}.
$$
(see, for example, \cite{Lnact}).
In this case, $D$ is identified with a countable dense subgroup 
of $\R.$  Note that $A$ has a unique tracial state $\tau$ such that
$\rho_A(K_0(A))=D.$
By the classification theorem (\cite{Lnmsri}), there is only one such
separable simple amenable \CA\, with tracial rank zero satisfying
the UCT.

Suppose that $\gamma_{0}: K_0(C(X))\to K_0(C(X))$ is an order
isomorphism and $\gamma_1: K_1(C(X))\to K_1(C(X))$ is an
isomorphism. Then $\gamma_0$ induces an order isomorphism
$I(\gamma_0): D\oplus H_{00}\to D\oplus  H_{00}$ such that
$I(\gamma_0)((d,0))=(d,0)$ and if $\kappa_{n,\infty}^{(0)}:
G_{00}\to H_{00}$ is the \hm\, induced by the inductive system,
$I(\gamma_0)\circ\kappa_{n,\infty}^{(0)}(g)=\kappa_{n,\infty}^{(0)}(\gamma_0(g))$
for all $g$ in the $n$-th $\Z\oplus G_{00}.$ Similarly
$I(\gamma_1)\circ\kappa_{n,\infty}^{(1)}(g)=\kappa_{n,\infty}\circ
\gamma_1(g)$ for all $g$ in the $n$-th $G_1.$
}

\end{Dfn}

\begin{thm}\label{N2MT}
Let $X$ be a connected compact metric space with
$K_0(C(X))=\Z\oplus G_{00}$ and $K_1(C(X))=G_1.$ Fix $\{a_n\}$ and
$\{b_n\}$ two sequences of positive numbers. Let $A=A(X, \{a_n\},
\{b_n\}).$ Suppose that $X$ has a minimal homeomorphism $\psi$
such that $(X,\psi)$ has mean dimension zero. Then there exists an
automorphism $\af: A\to A$ with the tracial cyclic Rokhlin
property such that
$$
\af_{*i}=I(\psi^{\natural}_{*i})\,\,\,{\rm on}\,\,\,K_i(A),
$$
$i=0,1.$
\end{thm}

\begin{proof}
Let $k(1)=1,$ $k(n+1) = b_nk(n),$ $n=1,2,...,.$ Put $l_n=b_n-a_n.$ We
may assume, as mentioned in Definition \ref{Construction},  that
 that $\lim_{n\to\infty}{a_n\over{b_n}}=0.$

Choose $\ep_n={1\over{2^n}}.$ Fix a dense  subset
$\{f_1,f_2,...,f_n,...\}.$ Let $\{\dt_n\}$ be as in Lemma \ref{Nlem1.2}
associated with $\{k(n)\},$ $\{f_1,f_2,...,f_n,...\}$ and
$\{\ep_n\}.$
 Let $F_n=\{x(1,n),x(2,n),...,x(l(n),n)\}$ be a
finite subset (with $l(n)$ elements) of $X$ such that $F_n$  is
$1/2^n$ -dense in $X$ and there is a permutation $s_n:
(1,2,...,l(n)) \to (1,2,...,l(n))$ such that
$$
{\rm dist}(x(j,n), \psi(x_{s_n(x(j,n))}))<\min\{\dt_n,\ep_n\}\,
j=1,2,...,l(n).
$$
Such $F_n$ is given by Lemma \ref{EVEN}.

Define $\phi_n: C(X)\otimes M_{k(n)}\to C(X)\otimes M_{k(n+1)}$ by
\begin{eqnarray}\label{urem}
\phi_n(f)={\rm diag}({\overbrace{f,f,...,f}^{a_n}},
f(x(1,n)),f(x(2,n)),...,f(x(l(n),n)))
\end{eqnarray}
for $f\in C(X)\otimes M_{k(n)}.$ Then $A\cong
\lim_{n\to\infty}(C(X)\otimes M_{k(n)},\phi_n).$

Note that $A$ has a unique tracial state (see Definition \ref{Construction}).
Denote it by $\tau.$


Define $\psi^{\natural} : C(X) \rightarrow C(X)$ by
$\psi^{\natural}(f)(x) = f(\psi(x))$ for $x \in X.$ Put $A_n=C(X)\otimes M_{k(n)}$ for each $n \in \N$ and 
define $\alpha_n = \psi^{\natural}\otimes id_{k_n}: A_n\to A_n.$  Then $\alpha_n
\in \Aut(A_n)$ for each $n \in \N$. We set
$$
{\mathcal F}_n = \{(a_{i,j})\in M_{k(n)}: a_{i,j}\in
\{f_1,f_2,...,f_n\}\cup {\mathbb C}\}.
$$
 From Lemma~$\ref{Nlem1.2}$ we may assume that
there is a permutation matrix $U_{n+1}' \in M_{l(n)}$ such that
\begin{eqnarray}\label{t2e1}
||\phi_n \circ \alpha_n(f) - {\rm ad}(U_{n+1})\circ
\alpha_{n+1}\circ \phi_n(f)|| < \ep_{n}
\end{eqnarray}
for each $n \geq 1$ and $f \in {\cal F}_n\cup
\cup_{j=1}^{n-1}\phi_j({\mathcal F}_j),$ where $U_{n+1}={\rm
diag}(1_{a_nk(n)}, 1_{k(n)}\otimes U_{n+1}')\in M_{k(n+1)}.$
Thus we obtain a permutation matrix $V_n'\in M_{l(n+1)}$ such that
the following approximately intertwining:

\begin{align*}
\begin{array}{ccccccc}
A_{1} & \mapright{\phi_{1}} &A_{2}&\mapright{\phi_{2}} &A_{3}
& \mapright{\phi_{3}}&\cdots \\
\mapdown{\alpha_1} &&\mapdown{{\rm ad}V_2\circ \alpha_2}
&&\mapdown{{\rm ad}V_3\circ \alpha_3}
&\\
A_{1} & \mapright{\phi_{1}} &A_{2}&\mapright{\phi_{2}} &A_{3}
& \mapright{\phi_{3}}&\cdots, \\
\end{array}
\end{align*}
where $V_n$ is a sequence of permutation matrixes in $M_{k_n}$
such that

\begin{eqnarray}\label{T2e2}
||\phi_n \circ {\rm ad}(V_n)\circ\alpha_n(f) - {\rm
ad}(V_{n+1})\circ \alpha_{n+1}\circ \phi_n(f)|| < \ep_{n},
\end{eqnarray}
for $f \in {\mathcal F}_n\cup \cup_{j=1}^{n-1}\phi_j({\cal F}_j)$.

Then there is a *-homomorphism $\alpha\colon A \rightarrow A$ such that
\begin{eqnarray}\label{T2e3}
\alpha(\phi_{l,\infty}(x)) =
\lim_{i\rightarrow\infty}
\phi_{i, \infty}\circ {\rm ad}(V_{i})\circ \alpha_i\circ \phi_{l,i}(x),
\end{eqnarray}
$x \in A_l$ and $l = 1,2,\dots$, where $\phi_{l,i} =
\phi_{i-1}\circ\phi_{i-2}\circ\cdots\circ\phi_l$. Moreover $\af$
is an automorphism.

It is also easy to check that
\begin{eqnarray}\label{T2eK}
\af_{*i}=I(\psi^{\natural})_{*i},\,\,i=0,1.
\end{eqnarray}

We also note that, by (\ref{T2e2}),
\begin{eqnarray}\label{T2e4}
\af(\phi_{m,\infty}(f))\approx_{1/2^{m-1}} \phi_{m,\infty}\circ
{\rm ad}\,V_m\circ \af_m(f)\,\,\,{\rm for}\,\,\, f\in {\cal F}_m.
\end{eqnarray}

Denote by $j_n: C(X)\to C(X)\otimes M_{k(n)}$ the embdedding
defined by $j_n(f)=f\otimes 1.$

Claim: for each $m\ge 1$ and $\tau\circ \phi_{m,\infty}\circ
j_m(f)=\tau\circ\phi_{m,\infty}\circ j_m(f\circ\psi)$ for all
$f\in C(X).$

Fix $m\ge 1.$ If $n\ge m+1,$ there is a projection $q_n\in
M_{k(n)}$ such that
$$
\phi_{m,n}\circ j_m(f) q_n= q_n\phi_{m,n}\circ
j_m(f),\,\,\,\,\,tr(q_n)={a_n\over{b_n}}
$$
and, one may write that
$$
\phi_{m,n}\circ j_m(f)(1-q_n)={\rm diag}
(f(x(1,n)),f(x(2,n),...,f(x(l(n),n)))
$$
for all $f\in C(X).$ Therefore
$$
\phi_{m,n}\circ j_m(f\circ \psi)(1-q_n)={\rm
diag}(f(\psi(x(1,n))),f(\psi(x(2,n))),..., f(\psi(x(l(n),n))))
$$
for all $f\in C(X).$ Denote by ${\bar q}_n=\phi_{n,\infty}(q_n).$
For any $\ep>0,$ and any $f\in C(X),$ by the virtue of
Lemma \ref{Nlem1.2}, if $n$ is sufficiently large, there is a
permutation matrix $W_n\in M_{k(n)}$ such that
$$
\|W_n^*[\phi_{m,n}\circ j_m(f)(1-q_n)]W_n-\phi_{m,n}\circ
j_m(f\circ\psi)(1-q_n)\|<\ep.
$$
Put ${\bar W}_n=\phi_{n,\infty}(W_n).$ The above implies that
$$
\|{\bar W}_n^*\phi_{m,\infty}\circ j_m(f)(1-{\bar q}_n)
{\bar W}_n-\phi_{m,\infty}\circ j_m(f\circ\psi)(1-{\bar q}_n)\|<\ep.
$$
Thus we conclude that
$$
|\tau(\phi_{m,\infty}\circ j_m(f))-\tau(\phi_{m,\infty}\circ
j_m(f\circ \psi))| <{a_n\over{b_n}}+\ep
$$
for all sufficiently large $n.$ Since
$\lim_{n\to\infty}{a_n\over{b_n}}=0.$ This implies that
$$
\tau(\phi_{m,\infty}\circ j_m(f))=\tau(\phi_{m,\infty}\circ
j_m(f\circ \psi))
$$
for all $f\in C(X).$ This proves the claim.

So
\begin{eqnarray}\label{T2m}
\tau\circ\phi_{m,\infty}\circ j_m(f)=\int_X f d\mu
\end{eqnarray}
for all $f\in C(X)$ and for some $\psi$-invariant Borel
probability measure $\mu$.

Now we verify that $\af$ has tracial Rokhlin property. For this
end, we fix an integer $N\ge 1,$ $\ep>0$ and a finite subset
${\cal F}\subset A.$ For convenience, without loss of generality,
we may assume that there exists ${\cal F}'\subset C(X)\otimes
M_{k(m)}$ such that
$$
\phi_{m,\infty}({\cal F}')\supset {\cal
F}\cup\cup_{j=1}^N\af^{-j}({\cal F})
$$
We may further assume that ${\cal F}'$ is in the unit ball of
$C(X)\otimes M_{k(m)}$ for some integer $m\ge 1.$ We may also
assume that
$$
{\cal F}'\subset \{(a_{i,j}) : a_{i,j} \in {\cal G}, 1 \leq i, j
\leq k_{i_l}\},
$$
where ${\cal G}\subset C(X)$ is a finite subset in the unit ball.
Without loss of generality, we may further assume that
$$
{\cal G}\subset \{f_1,f_2,...,f_m\}.
$$

Choose $\eta>0$ such that, if $x,x'\in X$ and ${\rm
dist}(x,x')<\eta,$
\begin{eqnarray}\label{T2e5}
|f(x)-f(x')| < \frac{\ep}{4N^28(k(m))^2}
\end{eqnarray}
for all $f\in {\cal G}.$

It follows from Lemma \ref{Lrok} that there are mutually disjoint open
subsets $G_1,G_2,...,G_L\subset X$ such that

(i) ${\rm diam}(G_i)<\eta,$

(ii) $\psi^j(\cup_{i=1}^LG)\cap
\psi^l(\cup_{i=1}^LG_i))=\emptyset,$ if $j\not=l$ and
$j,l= 0,1,...,N - 1,$

(iii) $\mu(\partial(G_i))=0$ for all $\mu\in T_{\psi}$ and

(iv) $\mu(\cup_{j=0}^{N-1} \psi^j(\cup_{i=1}^LG_i))>1-\ep/16$ for all
$\mu\in T_{\psi}.$

There is, for each $i,$ another open subset $S_i$ such that the
closure of $S_i$ is in $G_i,$ and
\begin{eqnarray}\label{T2eTr}
\mu(S_i)>\mu(G_i)-\ep/4LN\andeqn \mu(\partial(S_i))=0
\end{eqnarray}
for all $\mu\in T_{\psi}.$

Let $g_i'\in C(X)$ such that $0\le g_i'(x)\le 1,$ $g_i'(x)=0$ if
$x\not\in G_i,$ $0<g_i'(x)$ if $x\in G_i$ and $g_i'(x)=1$ if $x\in
{\bar S}_i,$ $i=1,2,...,L.$ Define $a_i'\in C(X)$ such that $0\le
a_i'(x)\le 1,$ $a_i'(x)>0$ for all $x\in S_i$ and $a_i'(x)=0$ if
$x\not\in S_i,$ $i=1,2,...,L.$

Fix $\eta_0>0.$ Choose $m_1\ge m$ such that
$1/2^{m_1}<\min\{\ep,\eta_0\}/16N^4.$ We also assume that
\begin{eqnarray}\label{T2e10}
{\rm dist}(a_i',\{f_1,f_2,...,f_{m_1}\})<\ep/4N\andeqn
{\rm dist}(g_i',\{f_1,f_2,...,f_{m_1}\})<\ep/4N
\end{eqnarray}
$i=1,2,...,L.$

 We view $g_i'$ and $a_i'$ as elements in
$E_{11}M_{k(m_1)}E_{11},$ where $\{E_{i,j}\}$ is a system of
matrix unit for $M_{k(m_1)}.$ Denote by
$a_i''=\phi_{m_1,m_1+1}(a_i')$ $a_i=\phi_{m_1,\infty}(a_i'),$
$g_i''=\phi_{m_1,m_1+1}(g_i')$ and ${\bar
g}_i=\phi_{m_1,\infty}(g_i'),$ $i=1,2,...,L.$

Consider the hereditary \SCA\, $\overline{a_iAa_i}.$ Then, by the
claim and by (\ref{T2m}), there is a measure $\mu\in T_{\psi}$
such that
\begin{eqnarray}\label{T2e11}
\mu(S_i)=k(m_1)\sup\{\tau(b): 0\le b\le 1, b\in
\overline{a_iAa_i}\},\,\,\,i=1,2,...,L.
\end{eqnarray}
Since $A$ has real rank zero, one obtains a projection $p_i\in
\overline{a_iAa_i}$ such that
\begin{eqnarray}\label{T2e12}
\tau(p_i)>{\mu(S_i)\over{k(m_1)}}-\ep/2LNk(m_1),\,\,\,i=1,2,...,L
\end{eqnarray}

Put ${\bar E}_{i,j}=\phi_{m_1,\infty}(E_{i,j}).$
Note $p_i\le {\bar E_{1,1}}.$
With $1_A=\sum_{i=1}^{k(m_1)}{\bar E}_{i,i},$ write
$$
e_1={\rm
diag}(\overbrace{\sum_{i=1}^Lp_i,\sum_{i=1}^Lp_i,\dots,\sum_{i=1}^Lp_i}^{k(m_1)}).
$$

By the choice of $\eta$ (see also (\ref{T2e5}) and Lemma \ref{LDig}---see also 3.2 of \cite{Lnscd}), 
it follows that
\begin{eqnarray}\label{T2e13}
 \|e_1 \phi_{m,\infty}(f)-\phi_{m,\infty}(f)e_1\|<\ep/4N^2
\end{eqnarray}
for all $f\in {\cal F}'.$ By the definition of ${\cal F}',$ we
have
\begin{eqnarray}\label{T2e14}
\|\af^j (e_1)f-f\af^j(e_1)\|<\ep/4N^2
\end{eqnarray}
for $0 \le j \le N $ and all $f\in {\cal F}.$

Put $h_j=\phi_{m_1,\infty}({\rm
diag}(\overbrace{\sum_{i=1}^Lg_i'\circ
\psi^{j-1},\sum_{i=1}^Lg_i'\circ\psi^{j-1},...,
\sum_{i=1}^Lg_i'\circ\psi^{j-1}}^{k(m_1)})),$
$j=1,2,...,N.$ Note that
\begin{eqnarray}\label{T2e15}
h_jh_i=h_ih_j=0\,\,\,\,{\rm if}\,\,\,i\not=j,\,i,j=1,2,...,N.
\end{eqnarray}

We compute that
\begin{eqnarray}\label{T2e16}
{\rm ad}\,V_{m_1} \af_{m_1}(h_1)=h_2.
\end{eqnarray}

 Note that $e_1\le h_1.$ It follows from
(\ref{T2e16}) and (\ref{T2e4}) that
\begin{eqnarray}\label{T2w17}
\|\af^j(h_1)-h_{j+1}\|<N
(1/2^{m_1-1})<2N(1/2^{m_1})<\min\{\ep,\eta_0\}/8N^3,
\end{eqnarray}
$j=1,2,...,N - 1.$ Therefore, by (\ref{T2e15}),
\begin{eqnarray}\label{T2e18}
\|\af^j(e_1)\af^i(e_1)\|<\min\{\ep,\eta_0\}/4N^3,\,\,\,i\not=j,\,i,j=0,1,...,N.
\end{eqnarray}
By choosing a small $\eta_0$ (which depends only on $\ep$ and
$N$), it follows from \cite[Lemma 2.5.6]{Lntext} 
that there are mutually orthogonal
projections $e_1,e_2,...,e_N$ such that
\begin{eqnarray}\label{T2e19}
\|\af^j(e_1)-e_{j+1}\|<\ep/N,\,\,\,j=0,1,...,N-1.
\end{eqnarray}
It follows that
\begin{eqnarray}\label{T2eF1}
\|\af(e_i)-e_{i+1}\|<\ep,\,\,\,i=1,2,...,N - 1
\end{eqnarray}
By (\ref{T2e14}) and (\ref{T2e19}), we also have
\begin{eqnarray}\label{T2eF2}
\|e_if-fe_i\|<\ep,\,\,\,i=1,2,...,N
\end{eqnarray}
for all $f\in {\cal F}.$

 By (\ref{T2e12}),
\begin{align*}
\tau(e_1) &= k(m_1)\tau(\sum_{i=1}^Lp_i)
> \sum_{i=1}^L(\mu(S_i) - \frac{\ep}{2LN})\\
&= \sum_{i=1}^L\mu(S_i) - \frac{\ep}{2N}
> \sum_{i=1}^L(\mu(G_i) - \frac{\ep}{4LN}) - \frac{\ep}{2N}\\
&= \mu(\cup_{i=1}^LG_i)-\frac{3\ep}{4N}.
\end{align*}
By (\ref{T2eK}),
\begin{eqnarray}\label{T2e21}
\tau(\af^j(e_1))=\tau(e_1),\,\,\,j=1,2,...,N-1.
\end{eqnarray}

 Finally, by (iv) and (\ref{T2eTr}), we compute that
\begin{align*}
\tau(\sum_{i=1}^Ne_i) &= N\tau(e_1)
> N\mu(\cup_{i=1}^LG_i) - \frac{3\ep}{4}\\
&= \mu(\cup_{j=0}^{N-1}\psi^j(\cup_{i=1}^LG_i)) - \frac{3\ep}{4}\\
&> 1 - \frac{\ep}{16} - \frac{3\ep}{4}
> 1-\ep.
\end{align*}

Hence $\af$ has the tracial Rokhlin property.
By (\ref{T2eK}), $\af_{*0}|_D = {\rm id}_{D}.$  Since $\rho_A(D)=\rho_A(K_0(A)),$ 
it follows from \cite[Theorem 3.16]{Lndyn} that $\af$ has the tracial cyclic Rokhlin property.
\end{proof}

\begin{Remark}
{\rm
In the proof of Theorem \ref{N2MT}, one  could work on more general 
inductive limits. For example,  instead of 
considering a single summand $M_{k(n)}(C(X)),$ one can consider unital simple 
\CA s which are inductive limits of \CA s with form $\oplus_{i=1}^{l(n)}M_{r(i,n)}(C(X))$ with each 
partial maps having the same form as (\ref{urem}). 
By doing that, more general dimensional
groups $D$ and more complicated $K_0(A)$ can be obtained. However,  we choose this much simple construction  to shed some 
light on the relation of the non-commutative dynamical systems and 
commutative dynamical systems  which could be buried by
some more complicated combinatory argument. 

It should also be pointed out that the resulted crossed product \CA\, $A\rtimes_{\af}\Z$ 
obtained from Theorem \ref{N2MT} has tracial rank zero, 
by Theorem \ref{gpots}, and $A\rtimes_{\af}\Z$
is a unital simple AH-algebra.
}
\end{Remark}

%
%
%

\vskip 3mm

\begin{Remark}\label{Rro}
{\rm In the case that $X=S^1,$ let $\theta \in \R\backslash\Q.$
Consider the rotation $\phi\colon X \rightarrow X$ by $\phi(z) = (\exp{2\pi i\theta})z.$
Then $\phi$ is minimal and uniquely ergodic.
Theorem \ref{N2MT} gives automorphisms $\af$ on unital simple $A\T$-algebras
with tracial cyclic Rokhlin property. 
A modification of the proof of Theorem \ref{N2MT} will give
an automorphism with the tracial Rokhlin property on every unital simple $A\T$-algebra with real rank zero. But all of them are approximately inner.
However, one can use these approximately inner automorphisms with Rokhlin
property in the proof of Theorem \ref{MT1} instead of applying a deeper result 
of N. C. Phillips (\cite{Phn}).
 It should be noted that, in  this case the crossed product algebras
$A \rtimes_\alpha \Z$ are  again a unital simple $A\T$ algebra
of real rank zero by {\rm \cite[Corollary 3.6]{LO}} (and \cite{Lngpots}).
}
\end{Remark}

Not all connected finite CW complex $X$ have minimal homeomorphisms.
One of the interesting examples is Furstenberg transformation on $\T^k.$

\begin{Ex}\label{Exfurstenberg}
{\rm Let $\theta \in \R\backslash\Q$, $d \in \Z\backslash\{0\}$,
and let $g_j\colon \T^j\rightarrow \T$ be a continuous map,
$j=1,2,...,k-1.$ The Furstenberg transformation is defined by
$$
\phi(z_1,z_2,...,z_k)=(z_1e^{2\pi i\theta},
g_1(z_1)z_2,...,g_{k-1}((z_1,z_2,...,z_{k-1}))z_k)
$$
$|z_j|=1,$ $j=1,2,...,k.$
Furstenberg showed that
$\phi$ is minimal, if all $g_j$ are non-trivial, namely,
$g_j(z_1,\dots,z_j)=z_j^{d_j}e^{2\pi i f_j(z_1)}$ for some continuous function
$f_j: \T\to \R$ with $d_j\not=0,$ $j=1,2,...,k-1.$
Thus $\phi^{\natural}_{*1}\not={\rm id}_{K_1(C(\T^k))}.$
Therefore the automorphism $\af$ constructed from $\psi$ on $A_{\T^k}$
is not approximately inner since $\af_{*1}\not={\rm id}_{K_1(A_{\T^k})}.$
(See for example \cite{Rj}.)
}
\end{Ex}

\end{document}